\documentclass{article}
\usepackage{amsmath,amssymb,theorem,cite}

\newtheorem{theorem}{Theorem}[section]   

\newtheorem{proposition}[theorem]{Proposition}

\newtheorem{remark}[theorem]{Remark}
\numberwithin{equation}{section}
\numberwithin{theorem}{section}
\newcommand{\mc}[1]{{\mathcal #1}}
\newcommand{\bb}[1]{{\mathbb #1}}
\newcommand{\NLT}{\mc N\!\mc L\!\mc T}

\begin{document}

\title{Resonances and $O$-curves in Hamiltonian systems}

\author{Paolo Butt\`a\footnote{Dipartimento di Matematica, 
Universit\`a di Roma ``La Sapienza'', P.le Aldo Moro 2, 
00185 Roma, Italy. E-mail: {\tt butta@mat.uniroma1.it,
piero.negrini@uniroma1.it.} Fax: +39--06--44701007.}
\and Piero Negrini$^*$}

\maketitle

\begin{abstract} 
We investigate the problem of the existence of trajectories asymptotic
to elliptic equilibria of Hamiltonian systems in the presence of 
resonances.  

\end{abstract}


   \section{Introduction}
   \label{sec:1}
   
In the recent Congress held in Saint Petersbourg, dedicated to the 
50$^\text{th}$ anniversary of A.M. Lyapunov, the V.V. Kozlov conference
\cite{Koz1} has been devoted to the so-called Lyapunov's First Method,
applied in particular to the problem of the existence of $O^+$-curves
(resp.\ $O^-$-curves) that  are  integral curves asymptotic in the 
future (resp.\ in the past) to the equilibria of Lagrangian systems.  
This problem is obviously related to the problem of the inversion of 
Dirichlet-Lagrange Theorem and the papers \cite{Koz,Pal} contain the 
first important results  for analytic potential functions with
degenerate critical point. We refer to \cite{Ne} for a review on 
further researches in this subject.

In the framework of Hamiltonian systems, a large number of papers have
been devoted to the study of the relationship between instability (as
well stability) of equilibria and resonances.  We quote here just some
essential references \cite{Cab, Kha,Mar,Sok}: in all these papers the
instability was proved by constructing suitable Cetaev's functions.

We were therefore stimulated to study, by means of the First Method,
the existence of $O$-curves of Hamiltonian systems in the presence 
of resonances. The Hamiltonian function we consider is supposed to have
a non degenerate elliptic fixed point $F$ (see Section~\ref{sec:2}). 
The corresponding Birkhoff normal form turns out to be a perturbation
of an integrable Hamiltonian function. Moreover the integrable system admits two distinct straight lines $r_\pm$ which are $O^\pm$-curves. Following the russian terminology, these orbits are called ``rays of
the model system''. Their  existence is the starting point to 
build up the Cetaev function for the complete Hamiltonian system 
(resonances of order 3 or 4 are considered in \cite{Kha}). In 
Section~\ref{sec:3} we prove that the complete Hamiltonian system
admits $O$-curves, that have as asymptotic direction one of the two
rays of the model system. 

To prove this result we analyze the Hamiltonian system in a
neighborhood of each ray. Then, by using a suitable set of variables
we build up a new differential system with a non degenerate hyperbolic
equilibrium. According to the choice of the ray, the relationship between these two differential systems allows us to conclude that the
local stable  manifold of this hyperbolic equilibrium corresponds
either to a family of $O^+$-curves or to a family of $O^-$-curves 
of the Hamiltonian system.

Before concluding, let us recall that the  existence of $O$-curves 
was analyzed by Lyapunov in the case of {\it analytic} differential 
systems. Many results on this subject can be found in the book of Zubov\cite{Zu}.  In particular, in Chapter III, perturbations of homogeneous
polynomial differential systems having rays are studied and $O$-curves
are obtained as power series of the variable $t^{-\alpha}$, $\alpha$ being a  positive rational number depending on the degree of the  homogeneous system. Moreover, in the general case, the coefficient of
each term in the expansion is found to be a polynomial function of
the variable  $\log t$.     

The Hamiltonian functions we consider here are just $C^r$ functions
(the integer $r$ depending on the degree of the resonance) therefore
we need to apply general hyperbolic theory, i.e.\ fixed point theorems
in suitable functional spaces. 

In the case $C^\infty$ the method of expansion in power series could 
still be used, at least at a formal level. In fact the existence of 
formal series representing $O$-curves is effective: a deep result by
Kutsnesov \cite{Ku} allows to conclude that there exist true $O$-curves
having the formal series as asymptotic expansion.

In conclusion, we can repeat here the words one can read at the end of
the introduction to Chapter III of the book of Zubov : ``the basic 
ideas of almost all the results are in the works of Lyapunov''.    

   \section{Preliminaries}
   \label{sec:2}
   
We start by considering  a Hamiltonian function $H \in 
C^\infty(\Omega, \bb R)$, $\Omega$ a domain in $\bb R^{2n}$,
containing the critical point $x = 0$. Precisely, we assume that 
$x = 0$ is an elliptic critical point, and therefore we write
\begin{equation*}
H(x) = H_{[2]}(x) + \mc H(x),
\end{equation*}
where
\begin{equation*}
H_{[2]}(x) = \sum_{i =1}^{n}\omega_i (x_i^2 + x_{i +n}^2),
\qquad \mc H(x) = O(|x|^3).
\end{equation*}
We consider the case in which the quadratic form $H_{[2]}$ has non trivial zeros. Moreover we assume the following hypotheses.
\begin{itemize}
\item[$(H_1)$] There exists an integer $N \ge 3$ such that, for any   
$h =(h_1, \ldots, h_n) \in \bb Z_+^n \setminus \{0\}$, $|h| := 
\sum_{i =1}^{n} h_i\le N -1$,  
\begin{equation*}
\langle\omega, h\rangle \ne 0.
\end{equation*}
\item[$(H_2)$] There exists a vector $k$, $k \in 
\bb Z_+^n\setminus\{0\}$, $|k| = N$, such that
\begin{equation*}
\langle\omega,k\rangle = 0.
\end{equation*}
Moreover, if $k'\in \bb Z_+^{n}\setminus\{0, k\}$ is such that 
$\langle\omega,k'\rangle = 0$ then $|k'| > M:= 3N-1$.
\end{itemize}
In other words we assume a unique direction of resonance along $e:=\frac{k}{|k|}$, up to  order $3N$.

Under the previous hypotheses, by means of a symplectic change of variables $x = y + \Phi(y)$ in a neighborhood of $x =
0$, we obtain  the Birkhoff form of the Hamiltonian \cite{Sie},  
\begin{equation}
\label{eq:P5}
H^B(I,\theta) = \sum_{j=1}^{\frac 32 N} H^B_{[j]}(I) +  
H^B_R(I, \langle \theta,k\rangle) +  \mc R^B(I,\theta),
\end{equation}
where we also introduced the action-angle variables 
$I =(I_1, \ldots, I_n)$, $\theta =(\theta_1, \ldots, \theta_n)$, 
with
\begin{equation*}
\left\{\begin{array}{l} y_i = \sqrt{2 I_i}\cos \theta_i \\
y_{i + n} = \sqrt{2 I_i}\sin\theta_i \end{array}\right.
\qquad i = 1, \ldots, n.
\end{equation*}
In \eqref{eq:P5} we have:
\begin{itemize}
\item $H^B_{[j]}(I)$, $j \in \big\{1,\ldots,\frac 32 N\big\}$, 
is a {\it homogeneous polynomial} of degree $j$ in the variables 
$(I_1, \ldots, I_n)$. In particular 
\begin{equation*}
H^B_{[1]} = \langle \omega,I\rangle.
\end{equation*}
\item $H^B_R$, the resonant term, takes the form
\begin{equation*}
H^B_R(I,\langle \theta, k\rangle) =  \sum_{r = 0}^N
H_r (I, \langle \theta, k\rangle),
\end{equation*}
where, for $r \le N-1$,
\begin{equation}
\label{eq:P7.2}
H_r \,  =
\sum_{j_1 +\ldots +j_n = r} 
\Big\{h^{(1,r)}_{j_1,\ldots,j_n}
\cos\langle \theta, k\rangle + h^{(2, r)}_{j_1,\ldots,j_n} 
\sin\langle \theta, k\rangle \Big\}
\prod_{\alpha=1}^n I_\alpha^{\frac{ k_\alpha}{2}+ j_\alpha}
\end{equation}
and
\begin{eqnarray*}
H_N & = & \sum_{j_1 + \ldots + j_n = N} 
\Big\{h^{(1,N)}_{j_1,\ldots,j_n}\cos\langle \theta, k\rangle + 
h^{(2, N)}_{j_1,\ldots,j_n} \sin\langle \theta, k\rangle 
\\ && + \, h^{(3,N)}_{j_1,\ldots,j_n} \cos 2\langle \theta, 
k\rangle  + h^{(4, N)}_{j_1,\ldots,j_n} \sin 2\langle \theta, 
k\rangle \Big\} \prod_{\alpha=1}^n 
I_\alpha^{\frac{ k_\alpha}{2}+ j_\alpha}. ~~~~~~~~~
\end{eqnarray*}
\item $\mc R^B(I, \theta)$ collects the higher order terms:
\begin{equation*}
\mc R^B(I,\theta) = O\big(|I|^{\frac{ 3 N +1 }{2}}\big).
\end{equation*}
\end{itemize}

We now make two further assumptions on the system.
\begin{itemize}
\item[$(H_3)$] $I = k$ is a {\it channel of instability}, that is
\begin{equation*}
H^B_{[j]}(k) = 0\quad \forall\, j \in 
\Big\{1,\ldots,\frac{N- 1}{2}\Big\}.
\end{equation*}
\item[$(H_4)$] Let $\Psi :\bb R\to \bb R$ be defined by setting
\begin{equation}
\label{eq:P10.1}
\Psi(\sigma) :=  H^B_{[\frac 12  N]}(k) + H_0(k, k_1^2\sigma), 
\quad \sigma\in\bb R,
\end{equation}
with $H^B_{[\frac 12 N]}$ and $H_0$ as in \eqref{eq:P5} and 
\eqref{eq:P7.2} respectively. We assume that $\Psi$ has a simple  
zero at  $\sigma = c.$
 \end{itemize}
\begin{remark}\rm
\label{rem:ar}
The assumption $(H_3)$ is satisfied if $N\le 4$, otherwise it implies 
that Arnol'd condition of non degeneracy \cite{Ar} is violated up 
to the order $N$.
\end{remark}

\begin{remark}\rm
\label{rem:ar1}
The general form of \eqref{eq:P10.1} can be written
\begin{equation*}
\Psi(\sigma) = A\cos[k_1^2(\sigma + \sigma_0)]  +  B.
\end{equation*}
The hypothesis $(H_4)$ implies that there exist two roots 
$c_{\pm}$ of $\Psi$ such that $\Psi'(c_{+})\Psi'(c_{-}) < 0$.  Let 
$\Psi'(c_{+})> 0$. Then the ray $r_+$ (resp.\ $r_-$) corresponds to 
$c_+$ (resp.\ $c_-$).
Obviously, if $N$ is odd then $B = 0$ and therefore $\Psi''(c_\pm) = 
-  k_1^4 \Psi(c_\pm) = 0$. This fact will be exploited in 
Section~\ref{sec:3}.
\end{remark}

We are now able to state the main result of the paper.

\begin{theorem}
\label{theo:1}
Under Hypotheses $(H_1)$, $(H_2)$, $(H_3)$, and $(H_4)$ the Hamiltonian
system admits a $(n-1)$-parameters family of $O^+$-curves  as well 
as a $(n-1)$-parameters family of $O^-$-curves.
\end{theorem}

The proof of the theorem is given in the next section. 

\begin{remark}\rm
\label{rem:regularity}
The theorem  still hold  if we assume $H \in C^r(\Omega, \bb R)$, 
$r=3 N + 1$.
\end{remark}

   \section{Proof of Theorem \ref{theo:1}}
   \label{sec:3}

Without loss of generality we may assume $k_1>0$. We then introduce 
the symplectic linear change of coordinates 
$(I,\theta) \to  (J,\psi) = (A_1 I, A_2 \theta)$ given by
\begin{eqnarray*}
J_1 = k_1 I_1,& \qquad J_\alpha  = k_1 I_\alpha - k_\alpha I_1, &
\quad \alpha = 2, \cdots, n,
\\ \psi_1 = \frac{\langle \theta, k\rangle }{k_1^2}, & \qquad
{\displaystyle \psi_\alpha = \frac{\theta_\alpha}{k_1}}, & 
\quad \alpha = 2, \cdots, n.
\end{eqnarray*} 
The Hamiltonian function $K(J, \psi) = H(A_1^{-1}J, A_2^{-1}\psi)$ 
then reads
\begin{equation*}
K(J, \psi) =  K(J) +  K^R(J, \psi_1)+  \mc R (J, \psi),
\end{equation*}
where
\begin{eqnarray}
\label{eq:obv}
&& K(J) = \sum_{j=1}^{\frac 32 N} K_{[j]}(J), 
\qquad K_{[j]}(J) = H^B_{[j]}(A_1^{-1}J),
\nonumber\\ && K^R(J,\psi_1) = \sum_{r = 0}^N
K_r (J, \psi_1), \qquad K_r (J, \psi_1) = 
H_r(A_1^{-1}J, k_1^2 \psi_1), 
\nonumber\\
&& \mc R (J,\psi) =  \mc R^B(A_1^{-1}J, A_2^{-1}\psi).
\end{eqnarray}
We next introduce the following notation for $n$-components vectors, 
\begin{equation*}
x = (x_1,\hat x), \qquad \hat x \equiv (x_2, \ldots, x_n).
\end{equation*}
The Hamiltonian system then reads
\begin{equation}
\label{eq:proof1}
\left\{\begin{array}{l}  
\dot J_1 \, = \, {\displaystyle -\frac{\partial K^{R}(J_1, \hat J, \psi_1)}{\partial \psi_1} -\frac{\partial \mc R(J_1, \hat J,  \psi_1, 
\hat \psi)}{\partial \psi_1}}, \\
\\ \dot {\hat J} \, = \, {\displaystyle- \frac{\partial 
\mc R(J_1, \hat J,  \psi_1, \hat \psi)}{\partial \hat\psi}}, 
\\ \dot \psi_1 \, = \, {\displaystyle\Omega_1(J_1, \hat J) +
\frac{\partial K^{R}(J_1,\hat J, \psi_1)}{\partial J_1} +
\frac{\partial \mc R(J_1, \hat J,  \psi_1, \hat \psi)}{\partial J_1}}, 
\\ \\ \dot {\hat \psi} \, = \, {\displaystyle \hat \Omega(J_1, \hat J)  
+\frac{\partial K^R(J_1,\hat J,  \psi_1)}{\partial \hat J} 
+\frac{\partial \mc R(J_1, \hat J,  \psi_1, \hat \psi)}{\partial 
\hat J }},\end{array}\right.
\end{equation}
where 
\begin{equation*}
\Omega_1(J) = \frac{\partial K(J)}{\partial J_1}, \qquad 
\hat \Omega(J) = \frac{\partial K(J)}{\partial \hat J}.
\end{equation*}
We observe that, by virtue of Hypotheses $(H_3)$ and $(H_4)$,
since $A_1^{-1}(J_1,0) = J_1 k_1^{-2}k$,
\begin{equation}
\label{eq:proof3}
\Omega_1(J_1,0) + \frac{\partial K^R(J_1, 0, c)}{\partial J_1} =
\sum_{j=\frac N2 +\delta}^{\frac 32 N} 
\frac{\partial K_{[j]}}{\partial J_1}(J_1,0) + 
\sum_{r = 1}^N \frac{\partial K_r (J_1,0,c)}{\partial J_1},
\end{equation}
where 
\begin{equation*}
\delta = \left\{\begin{array}{ll} 1/2 & 
\text{ if $N$ is odd,} \\ 1 & \text{ if $N$ is even,}
\end{array}\right.
\end{equation*}
and
\begin{equation}
\label{eq:proof7}
\gamma := \frac{\partial K_0}{\partial \psi_1}(1,0, c)\ne 0. 
\end{equation} 
For definiteness we consider the case $\gamma > 0$; the case 
$\gamma <0$ will be briefly discussed at the end of the section. 
We look for solutions to (\ref{eq:proof1}) of the form
\begin{eqnarray}
\label{eq:n1}
J_1(t) &=& \left[\frac{2}{(N-2)\,t}\right]^{\frac{2}{N-2}}
\left\{\Gamma + u_1(t)\right\},
\nonumber\\
\hat J(t) &=&\left[\frac{2}{(N-2)\,t}\right]^{\frac{N +2}{N -2}}
\hat u(t),
\nonumber\\
\psi_1(t)&=& c  + \left[\frac{2}{(N-2)\,t}\right]^{
\frac{2\delta}{N-2}} \left\{c_1 + c_2 
\left[\frac{2}{(N-2)\,t}\right]^{\frac{2\delta}{N-2}} + v_1(t) 
\right\}, \nonumber\\
\hat \psi(t)&=& \left[\frac{2}{(N-2)\,t}\right]^{-1}
\big\{\hat \Omega_0 + \hat v(t)\big\},
\end{eqnarray}
where 
\begin{equation*}
\Gamma:=\gamma^{-\frac{2}{N -2}}, \qquad 
\hat \Omega_0 := \frac{2}{k_1(N -2)}(\omega_2, \ldots, \omega_n)
\end{equation*}
and $u(t) = (u_1(t),\hat u(t))$, $v(t) = (v_1(t),\hat v(t))$ are such 
that
\begin{equation*}
\lim_{t\rightarrow \infty}u(t)=0, \qquad
\lim_{t\rightarrow \infty}v(t)=0.
\end{equation*}
The explicit values of the parameters $c_1$ and $c_2$ are given in Proposition 3.1 below.  We introduce the new 
independent variable $z$ by setting
 \begin{equation*}
z = \left[\frac{2}{(N-2)\,t}\right]^{\frac{2}{N-2}}
 \end{equation*}
and write the differential equations for the functions
\begin{equation}
\label{eq:lem3.1}
\xi(z) := u(t(z))), \qquad\eta(z) := v(t(z)).
\end{equation}
This is the content of the following proposition.

\begin{proposition}
\label{prop:s}
Let, for $N$ odd,
\begin{equation}
\label{c1c2odd}
c_1 = -\frac{2\,\Gamma^{\frac{N-1}2}}{N+1}\frac{\partial K_{[\frac{N+1}2]}(1,0)}{\partial J_1},
\qquad c_2 = -\frac{2\,\Gamma^{\frac N2}}{N+2}
\frac{\partial K_1(1,0,c)}{\partial J_1} \end{equation}
and, for $N$ even,
\begin{equation}\label{c1c2even}
c_1 = -\frac{2\,\Gamma^{\frac N2}}{N+2} \left\{
\frac{\partial K_{[\frac N2 +1]}(1,0)}{\partial J_1}
+ \frac{\partial K_1(1,0,c)}{\partial J_1}\right\},
\qquad c_2 = 0.
\end{equation}
The functions $( \xi(z), \eta(z))$ are solutions to the differential system
\begin{equation}
\label{eq:lem5}
\left\{\begin{array}{l}
{\displaystyle z\frac{d\xi_1}{d z}= \frac{N - 2}{2} \xi_1 + d_0 z
+ U_1(z, \xi, \eta)}, \\ \\
{\displaystyle z \frac{d\hat \xi}{d z}= -\frac{N + 2}{2}\hat \xi + 
\hat U(z, \xi, \eta)}, \\ \\ 
{\displaystyle z \frac{d\eta_1}{d z} = -\frac{N + 2\delta}{2}
\eta_1 + d_1 z + d_2 \xi_1 + V_1(z, \xi, \eta)}, \\ \\
{\displaystyle z \frac{d \hat \eta}{d z} = 
\frac{N - 2}{2}\hat \eta + \hat d z + \hat V(z, \xi, \eta)},
\end{array}\right.
\end{equation}
where, for $N$ odd,
\begin{eqnarray}
\label{eq:lem6}
d_0&=& \Gamma^{\frac{N + 2}2}\, 
\frac{\partial K_1(1, 0, c)}{\partial \psi_1} - \Gamma \frac{k_1^4}2 c_1^2, 
\nonumber\\ 
d_1&=& - \Gamma^{\frac{N+1}2}
\frac{\partial K_{[\frac{N+3}2]}(1, 0) }{\partial J_1} -
c_1 \Gamma^{\frac N2} 
\frac{\partial^2 K_1(1, 0, c)}{\partial J_1\partial\psi_1},
\nonumber\\ 
d_2 &=& - \frac{N-1}2 \Gamma^{\frac{N-3}2}\frac{\partial 
K_{[\frac{N +1}2]}(1,0)}{\partial J_1} - c_1
\Gamma^{-1}\frac{N(N-2)}4, 
\end{eqnarray}
and, for $N$ even,
\begin{eqnarray}
\label{eq:lem61}
d_0&=& \Gamma^{\frac{N + 2}2}\, 
\frac{\partial K_1(1, 0, c)}{\partial \psi_1}+ c_1
\frac{\partial^2 K_0(1,0,c)}{\partial\psi_1^2}, 
\nonumber\\ 
d_1&=& - \Gamma^{\frac N2 +1} \left\{
\frac{\partial K_{[\frac N2+2]}(1, 0) }{\partial J_1} +
\frac{\partial K_2(1,0,c)}{\partial J_1} \right\} -
c_1 \Gamma^{\frac N2} 
\frac{\partial^2 K_1(1, 0, c)}{\partial J_1\partial\psi_1}
\nonumber \\ && + \,
k_1^4 c_1^2\, \Gamma^{\frac{N-2}2} 
\frac{\partial K_0(1,0,c)}{\partial J_1},
\nonumber\\ 
d_2 &=& - \frac N2 \Gamma^{\frac N2 -1}\frac{\partial 
K_{[\frac N2 +1]}(1,0)}{\partial J_1} - c_1
\Gamma^{-1}\frac{N(N-2)}4. 
\end{eqnarray}
Finally, for any $N$,
\begin{equation}
\label{eq:lem62}
\hat d = \Gamma^2 \left\{ 
\frac{\partial K_{[2]}(1, 0)}{\partial \hat J} + \delta_{N, 4}
\frac{\partial K_0(1, 0, c)}{\partial \hat J} \right\}.
\end{equation}
In particular $\hat d = 0$ if $N= 3$.

The functions $U=(U_1, \hat U)$ and $V= (V_1,\hat V)$ are 
$C^1$ functions in a  ``right neighborhood''
$\mc N^+$ of zero in $\bb R^{2n +1}$,
$$ 
\mc N^+ := \big\{(z, \xi, \eta) \in \bb R^{2n +1}\, :  \:
z \in [0, \varepsilon), \, \|(\xi, \eta)\| < \varepsilon \big\}.
$$
Moreover the Jacobian matrices satisfy
\begin{equation*}
DU(0,0,0) = 0, \qquad DV(0, 0, 0) = 0.
\end{equation*}
\end{proposition}

\noindent{\bf Proof.}
The proof is straightforward but cumbersome. Of course one has to rewrite the system \eqref{eq:proof1} in terms of the independent 
variable $z$ and the unknowns $(\xi, \eta)$, according to 
\eqref{eq:lem3.1}. The main point is then to extract from the r.h.s.\ 
of the system  the linear part and verify the regularity property 
of the remainder. 

By \eqref{eq:n1} and \eqref{eq:lem3.1} we have
\begin{eqnarray}
\label{tras}
&& J_1 = z(\Gamma+\xi_1),\qquad \hat J = z^{\frac{N+2}2} \hat \xi, \qquad
\psi_1 = c + z^\delta\big(c_1 + c_2 z^\delta + \eta_1\big),
\nonumber \\&&  \hat \psi = z^{-\frac{N-2}2}
\big(\hat \Omega_0 +\hat \eta \big),
\end{eqnarray}
whence, since $\dot z = - z^{\frac N2}$,
\begin{equation}
\label{eq:lem7}
\left\{\begin{array}{l}
{\displaystyle - z^{-\frac N2} \dot J_1 = \Gamma+\xi_1 + 
z\frac{d\xi_1}{d z}}, \\ \\
{\displaystyle - z^{-\frac N2} \dot {\hat J} = \frac{N+2}2 
z^{\frac N2} \hat\xi + z^{\frac{N+2}2} \frac{d\hat \xi}{d z}}, \\ \\ 
{\displaystyle - z^{-\frac N2} \dot \psi_1 = c_1\delta z^{\delta-1}
+ 2c_2 \delta z^{2\delta-1} + \delta z^{\delta-1} \eta_1 + 
z^\delta \frac{d\eta_1}{d z}}, \\ \\
{\displaystyle - \dot{\hat \psi} = - \frac{N-2}2
\big(\hat \Omega_0 +\hat \eta \big) + z \frac{d \hat \eta}{d z}}.
\end{array}\right.
\end{equation}
By comparing \eqref{eq:lem5} and \eqref{eq:lem7} it follows that
\begin{eqnarray}
\label{U=}
&& U_1(z, \xi, \eta) = -z^{-\frac N2} \dot J_1 - \Gamma - 
\frac N2 \xi_1 - d_0 z, \qquad
 \hat U(z, \xi, \eta) = - z^{-N} \dot{\hat J}, \quad\qquad \\ 
\label{V1=} && V_1(z, \xi, \eta) = -z^{1-\delta-\frac N2}\dot\psi_1 -
c_1\delta - 2c_2\delta z^\delta + \frac N2 \eta_1 - d_1 z - d_2\xi_1 ,  
\\ \label{Vh=} && \hat V(z, \xi, \eta) = 
-\dot{\hat\psi} + \frac{N-2}2 \hat\Omega_0 -\hat d z, \quad\qquad
\end{eqnarray}
where the time derivatives $(\dot J,\dot \psi)$ are given by
the r.h.s.\ of \eqref{eq:proof1} expressed in terms of the variables 
$(z,\xi,\eta)$ by means of \eqref{tras}.

In the sequel we shall denote by $\NLT$ a generic $C^1$ function of 
$(z,\xi,\eta)$ which vanishes with its first partial derivatives in
$(0,0,0)$. Since
\begin{equation}
\label{pp}
A_1^{-1}J = \frac{1}{k_1} J + \frac{J_1}{k_1^2}(0,\hat k)
= z \left\{\frac{\Gamma+\xi_1}{k_1^2} \, (k_1,\hat k) + z^{\frac N2} 
\frac{1}{k_1}(0,\hat\xi) \right\},
\end{equation}
recalling the definition of $\mc R$ in \eqref{eq:obv} we have:
\begin{equation}
\label{eq:pieroin1}
\frac{\partial \mc R(z(\Gamma+\xi_1),z^{\frac{N+2}2}\hat\xi, 
c+ \sqrt z (c_1 + c_2\sqrt z + \eta_1), z^{-\frac{N-2}2}
(\hat \Omega_0 +\hat \eta))}{\partial\psi} = z^{ N} \,\NLT, 
\end{equation}
\begin{equation}
\label{eq:pieroin2}
\frac{\partial \mc R(z(\Gamma+\xi_1),z^{\frac{N+2}2}\hat\xi, 
c+ \sqrt z (c_1 + c_2\sqrt z + \eta_1), z^{-\frac{N-2}2}
(\hat \Omega_0 +\hat \eta))}{\partial J} = z^{\frac N2} \,\NLT. 
\end{equation}
By the second equality in \eqref{U=} and \eqref{eq:pieroin1} we get
$\hat U= \NLT$. Moreover, by \eqref{Vh=}, \eqref{eq:pieroin2} and 
the definition \eqref{eq:lem62} of $\hat d$ it is straightforward to
conclude also that $\hat V = \NLT$, we omit the details. The analysis
of the functions $U_1$ and $V_1$ is more delicate and the cases 
$N$ odd and $N$ even have to be treated separately.

\smallskip
\noindent {\it Case $N$ odd.} 
By Remark~\ref{rem:ar1} we have
$$
\frac{\partial^2 K_0(1,0, c)}{\partial\psi_1^2} = 0, \qquad
\frac{\partial^3 K_0(1,0, c)}{\partial\psi_1^3} = 
- k_1^4\frac{\partial K_0(1,0, c)}{\partial\psi_1}.
$$
By the definitions \eqref{eq:obv}, \eqref{pp}, and recalling 
$\delta=\frac 12$ in this case, 
\begin{eqnarray*}
&& \frac{\partial K_0(z(\Gamma+\xi_1),z^{\frac{N+2}2}\hat\xi, c+ 
\sqrt z (c_1 + c_2\sqrt z + \eta_1))}{\partial\psi_1} \\ && 
~~~~~~~~~~~~~~~~~~~~~~~~~~~~~~ = 
z^{\frac N2} \frac{\partial K_0(\Gamma+\xi_1,z^{\frac N2}
\hat \xi, c + \sqrt z (c_1 + c_2\sqrt z + \eta_1))}{\partial\psi_1} 
\\ && 
~~~~~~~~~~~~~~~~~~~~~~~~~~~~~~ = 
z^{\frac N2} (\Gamma+\xi_1)^{\frac N2}
\frac{\partial K_0(1,0, c)}{\partial\psi_1}
\left\{ 1 - \frac{k_1^4}2 c_1^2 z + \NLT \right\},
\end{eqnarray*}
while, for $r\ge  1$,
\begin{eqnarray*}
&& \frac{\partial K_r(z(\Gamma+\xi_1),z^{\frac{N+2}2}\hat\xi, c+ 
\sqrt z (c_1 + c_2\sqrt z + \eta_1))}{\partial\psi_1}
\\ &&
~~~~~~~~~~~~~~~~~~~~~~~~~~~~~~ = 
z^{\frac N2 + r}\frac{\partial K_r(\Gamma+\xi_1,z^{\frac N2}
\hat \xi, c + \sqrt z (c_1 + c_2\sqrt z + \eta_1))}{\partial\psi_1} 
\\ & & 
~~~~~~~~~~~~~~~~~~~~~~~~~~~~~~ =  
z^{\frac N2 + r} (\Gamma+\xi_1)^{\frac N2+r} 
\frac{\partial K_r(1,0, c)}{\partial\psi_1} +  z^{\frac N2} \NLT.
\end{eqnarray*}
Then, recalling the definitions \eqref{eq:proof7}, \eqref{eq:lem6},
and $\Gamma=\gamma^{-\frac{2}{N -2}}$,
\begin{eqnarray*}
-z^{-\frac N2} \dot J_1 & = & \gamma(\Gamma+\xi_1)^{\frac N2} 
\left(1 - \frac{k_1^4}2 c_1^2 z \right) + 
\frac{\partial K_1(1,0, c)}{\partial\psi_1} \Gamma^{\frac {N+2}2} z
+ \NLT \\
& = & \Gamma + \frac N2 \xi_1 + d_0 z + \NLT.
\end{eqnarray*}
The previous expansions imply 
$U_1 = \NLT$. 

Analogously, recalling also \eqref{eq:proof3},
\begin{eqnarray*}
\dot\psi_1 & = & z^{\frac{N-1}2} \bigg\{
\frac{\partial K_{[\frac{N+1}2]} (\Gamma+\xi_1,0)}{\partial J_1} + 
(c_1+c_2 \sqrt z + \eta_1) \frac{\partial^2 K_0(\Gamma+\xi_1,0,c)}
{\partial J_1\partial\psi_1} \bigg\}
\\ && + \, z^{\frac N2} 
\frac{\partial K_1(\Gamma+\xi_1,0,c)}{\partial J_1}  + 
c_1^2 z^{\frac{N+1}2}
\frac{\partial^3 K_0(\Gamma+\xi_1,0,c)}{\partial J_1\partial\psi_1^2}
\\ && \, + \, c_1 z^{\frac{N+1}2}
\frac{\partial^2 K_1(\Gamma+\xi_1,0,c)}{\partial J_1\partial\psi_1}
+ z^{\frac{N+1}2}
\frac{\partial K_{[\frac{N+3}2]}(\Gamma+\xi_1,0)}{\partial J_1} 
\\ && +\, z^{\frac {N-1}2} \NLT,
\end{eqnarray*}
whence
\begin{eqnarray*}
z^{\frac{1-N}2} \dot\psi_1 & = & 
(\Gamma+\xi_1)^{\frac{N-1}2} 
\frac{\partial K_{[\frac{N+1}2]}(1,0)}{\partial J_1} + 
c_1 (\Gamma+\xi_1)^{\frac{N-2}2}\frac{\partial^2 K_0(1,0,c)}
{\partial J_1\partial\psi_1} 
\\ && + \, (c_2\sqrt z + \eta_1) 
(\Gamma+\xi_1)^{\frac{N-2}2}\frac{\partial^2 K_0(1,0,c)}
{\partial J_1\partial\psi_1} 
\\ && + \, \sqrt z \, \Gamma^{\frac N2} 
\frac{\partial K_1(1,0,c)}{\partial J_1} + 
c_1^2 z\, \Gamma^{\frac{N-2}2}
\frac{\partial^3 K_0(1,0,c)}{\partial J_1\partial\psi_1^2}
\\ && + \, c_1 z\, \Gamma^{\frac N2}
\frac{\partial^2 K_1(1,0,c)}{\partial J_1\partial\psi_1}
+ z\, \Gamma^{\frac{N+1}2} 
\frac{\partial K_{[\frac{N+3}2]}(1,0)}{\partial J_1} 
\\ && +\, z^{\frac {N-1}2} \NLT. 
\end{eqnarray*}
On the other hand, by the explicit form of the functions $K_0$, 
recalling $|k| = N$ and Remark~\ref{rem:ar1},
it is easy to verify that
\begin{equation*}
\frac{\partial^2 K_0(1,0,c)}
{\partial J_1\partial\psi_1} = \frac N2 \frac{\partial K_0(1,0,c)}
{\partial\psi_1} = \frac N2 \gamma, \quad
\frac{\partial^3 K_0(1,0,c)}{\partial J_1\partial\psi_1^2}
= - k_1^4 \frac{\partial K_0(1,0,c)}{\partial J_1} = 0.
\end{equation*} 
By inserting the previous expression of $z^{\frac{1-N}2} \dot\psi_1$ 
in \eqref{V1=}, expanding up to the first order in the variable
$\xi_1$, and recalling the definitions \eqref{c1c2odd},
\eqref{eq:lem6} of $c_1$, $c_2$, $d_1$, and $d_2$,
we get $V_1=\NLT$. 

\smallskip
\noindent {\it Case $N$ even.} 
To prove $U_1=\NLT$ we argue as before; the only difference is
that $\delta=1$ and $c_2=0$ in this case, so that
\begin{eqnarray*}
&& \frac{\partial K_0(z(\Gamma+\xi_1),z^{\frac{N+2}2}\hat\xi, c+ 
z (c_1 + \eta_1))}{\partial\psi_1} \\ && 
~~~~~~~~~~~~~~~ = 
z^{\frac N2} (\Gamma+\xi_1)^{\frac N2}
\left\{\frac{\partial K_0(1,0, c)}{\partial\psi_1}+ 
\frac{\partial^2 K_0(1,0, c)}{\partial\psi_1^2} c_1 z 
+ \NLT \right\},
\end{eqnarray*}
whence the definition of $d_0$ in \eqref{eq:lem61} for $N$ even.

We finally have
\begin{eqnarray*}
\dot\psi_1 & = & z^{\frac N2} \bigg\{
\frac{\partial K_{[\frac N2+1]} (\Gamma+\xi_1,0)}{\partial J_1} + 
(c_1 + \eta_1) \frac{\partial^2 K_0(\Gamma+\xi_1,0,c)}
{\partial J_1\partial\psi_1} \bigg\}
\\ && + \, z^{\frac N2} 
\frac{\partial K_1(\Gamma+\xi_1,0,c)}{\partial J_1}  + 
+ c_1^2 z^{\frac N2+1}
\frac{\partial^3 K_0(\Gamma+\xi_1,0,c)}{\partial J_1\partial\psi_1^2}
\\ && +\, c_1 z^{\frac N2 +1}
\frac{\partial^2 K_1(\Gamma+\xi_1,0,c)}{\partial J_1\partial\psi_1}
+ z^{\frac N2 +1} \frac{\partial K_{[\frac N2 +2]}
(\Gamma+\xi_1,0)}{\partial J_1} 
\\ && + \, z^{\frac N2 +1} \frac{\partial K_2
(\Gamma+\xi_1,0,c)}{\partial J_1} + z^{\frac N2} \NLT,
\end{eqnarray*}
whence
\begin{eqnarray*}
z^{-\frac N2} \dot\psi_1 & = & 
(\Gamma+\xi_1)^{\frac N2} \frac{\partial K_{[\frac N2 +1]} (1,0)}{\partial J_1} + c_1 (\Gamma+\xi_1)^{\frac{N-2}2}
\frac{\partial^2 K_0(1,0,c)}{\partial J_1\partial\psi_1} 
\\ && + \, \Gamma^{\frac N2} \frac{\partial K_1(1,0,c)}{\partial J_1}
+ \eta_1 (\Gamma+\xi_1)^{\frac{N-2}2}
\frac{\partial^2 K_0(1,0,c)}{\partial J_1\partial\psi_1} 
\\ && + \, c_1 z\,\Gamma^{\frac N2}
\frac{\partial^2 K_1(1,0,c)}{\partial J_1\partial\psi_1}
+ c_1^2 z\, \Gamma^{\frac{N-2}2}
\frac{\partial^3 K_0(1,0,c)}{\partial J_1\partial\psi_1^2} 
\\ && +\, z\, \Gamma^{\frac N2 +1}
\frac{\partial K_{[\frac N2 +2]}(1,0)}{\partial J_1} + 
z\, \Gamma^{\frac N2 +1} \frac{\partial K_2(1,0,c)}{\partial J_1} + \NLT.
\end{eqnarray*}
By inserting the previous expression of $z^{-\frac N2}\dot\psi_1$ in 
\eqref{V1=}, expanding up to the first order in the variable
$\xi_1$, using that in this case we still have
\begin{equation*}
\frac{\partial^2 K_0(1,0,c)}
{\partial J_1\partial\psi_1} = \frac N2 \frac{\partial K_0(1,0,c)}
{\partial\psi_1} = \frac N2 \gamma, \quad
\frac{\partial^3 K_0(1,0,c)}{\partial J_1\partial\psi_1^2}
= - k_1^4 \frac{\partial K_0(1,0,c)}{\partial J_1},
\end{equation*} 
and recalling the definitions \eqref{c1c2even} and \eqref{eq:lem61} 
of $c_1$, $c_2$, $d_1$, and $d_2$, we get $V_1=\NLT$.
\hfill $\square$\medskip

System (\ref{eq:lem7}) is equivalent to the autonomous system:
\begin{equation}
\label{eq:lem5.1}
\left\{\begin{array}{l}
{\displaystyle \frac{dz}{d\tau} = - z}, \\ \\
{\displaystyle \frac{d\xi_1}{d \tau}= 
- \frac{N - 2}{2} \xi_1 - d_0 z - U_1(z, \xi, \eta)}, \\ \\
{\displaystyle \frac{d\hat \xi}{d \tau}= 
\frac{N + 2}{2}\hat \xi - \hat U(z, \xi, \eta)}, \\ \\ 
{\displaystyle \frac{d\eta_1}{d \tau} = 
\frac{N + 2\delta}{2}\eta_1 - d_1 z - d_2\xi_1- V_1(z, \xi, \eta)},
\\ \\ {\displaystyle \frac{d \hat \eta}{d \tau} = 
- \frac{N - 2}{2}\hat \eta - \hat d z - \hat V(z, \xi, \eta)}.
\end{array}\right.
\end{equation}
The origin  $(z, \xi, \eta) = (0, 0, 0)$ is an hyperbolic equilibrium.
In order to apply the hyperbolic theory we have to get rid of the fact
that the system is defined in $\mc N^+$, and not in a full neighborhood
of the origin. However the proof of the existence of the stable 
manifold can be easily adapted to this case.  Then we have a local
invariant stable  manifold, parametrized by $(z, \xi_1, \hat \eta)$, which is the graph of a $C^1$ function defined in a neighborhood 
$B_\varepsilon \subset \mc N^+$.  

In conclusion, coming back to the original variables we have obtained
$W^+$, an $(n-1)$-dimensional surface of $O^+$-curves.

We conclude this section by briefly considering the case $\gamma < 0$.
We consider the  problem \eqref{eq:proof1} in the past, or
equivalently, we change the sign of the Hamiltonian function and let 
$t > 0$. In particular we replace $\gamma$ by $|\gamma|$ everywhere it
appears. We thus get again a system of the same form of 
\eqref{eq:lem5.1}. In conclusion, to the  stable manifold corresponds
now $W^-$, an $(n-1)$-dimensional surface of  $O^-$-curves.

Finally, recalling that we have two different rays of the model
system, the proof of Theorem \ref{theo:1} is accomplished.

\end{document}